\documentclass{pja00}

\usepackage{amsthm, amssymb, amsfonts, mathrsfs, amsmath}
\usepackage[T5]{fontenc}

\usepackage{listings}

\theoremstyle{definition}
\newtheorem{defn}{Definition}[section]

\newtheorem*{acknow}{Acknowledgments}

\theoremstyle{plain}
\newtheorem{thm}[defn]{Theorem}

\newtheorem{corl}[defn]{Corollary}

 \runninghead{}

\usepackage[x11names]{xcolor}
\usepackage[
  colorlinks,
]{hyperref}

\AtBeginDocument{\hypersetup{citecolor=magenta, urlcolor=magenta, linkcolor=blue}}

\title{On the lambda algebra and Singer's cohomological transfer}

\Author{\DJ\d \abreve ng}{V\~o Ph\'uc}

\affiliation{1}{Faculty of Education Studies, University of Khanh Hoa, Vietnam}


\begin{document}

\maketitle
\begin{abstract}
Writing $\mathbb A$ for the 2-primary Steenrod algebra, which is the algebra of stable natural endomorphisms of the mod 2 cohomology functor on topological spaces. Working at the prime 2, computing the cohomology of $\mathbb A$ is an important problem of Algebraic topology,  because it is the initial page of the Adams spectral sequence converging to stable homotopy groups of the spheres. A relatively efficient tool to describe this cohomology is the Singer algebraic transfer of rank $n$ in \cite{Singer}, which passes from a certain subquotient of a divided power algebra to the cohomology of $\mathbb A.$ Singer predicted that this transfer is a monomorphism, but this remains open for $n\geq 4.$ This short note is to verify the conjecture in the ranks 4 and 5 and some generic degrees.

\end{abstract}

\section{Introduction}\label{s1} 
It is well-known that there is a group homomorphism $Sq^{n}$ for $n\geq 0,$ between mod-2 cohomology groups of a topological space, called Steenrod squares of degrees $n$. They are stable cohomology operations, that is, they commute with suspension maps. All the Steenrod squares form an algebraic structure which is known as the 2-primary Steenrod algebra $\mathbb A$ subject to the Adem relations. 
%Milnor \cite{J.M} also showed that $\mathbb A$ has a filtration by finite dimensional sub algebras and is a finite dimensional vector space in each grading and so, the dual algebra also has a Hopf algebraic structure. 
It was applied to the vector fields on spheres and the Hopf invariant one problem, which asks for which $n$ there exist maps of Hopf invariant $\pm 1$. 
%In addition, the Steenrod squares are also used for computing homotopy groups of spheres. 
So, the Steenrod algebra is one of the important tools in Algebraic topology. Specifically, its cohomology ${\rm Ext}_{\mathbb A}^{*, *}(\mathbb Z/2, \mathbb Z/2)$ is an algebraic object that serves as the input to the Adams spectral sequence (the ASS) \cite{J.A} and therefore, computing this cohomology is of fundamental importance to the study of the stable homotopy groups of spheres. 

The cohomological transfer, defined by Singer \cite{Singer}, could be an useful approach to describe the mysterious structure of the cohomology algebra of $\mathbb A$. In order to better understand this transfer, we will use the following notations and the relevant concepts. Let denote $V^{\oplus n}$ the $n$-dimensional vector space over the prime field $\mathbb Z/2.$ Then, we write $H^{*}(V^{\oplus n})$ and $H_*(V^{\oplus n})$ for mod-2 cohomology  and homology of $BV^{\oplus n}$ (the classifying space of $V^{\oplus n}$). One should note that $BV^{\oplus n}$ is homotopy equivalent to the cartesian product of $n$ copies of the union of the finite projective spaces. As it is known, $H^{*}(V^{\oplus n})$ is identified with the polynomial algebra $\mathbb Z/2[u_1, \ldots, u_n]$ on generators of degree $1,$ equipped with the canonical unstable algebra structure over the Steenrod algebra (i.e., it is a commutative, associative, graded $\mathbb Z/2$-algebra equipped with a structure of unstable $\mathbb A$-module and satisfying two relations, one called the Cartan formula, and the other called the instability relation: $Sq^{\deg(x)}(x) = x^{2}$.) By dualizing, $H_*(V^{\oplus n})$ has a natural basis dual to the monomial basis of $H^{*}(V^{\oplus n}).$ We denote by $a_1,\ldots, a_n$ the basis of $H_1(V^{\oplus n})$ dual to the basis $u_1, \ldots, u_n$ of $H^{1}(V^{\oplus n}) ={\rm Hom}(V^{\oplus n}, \mathbb Z/2),$ so that \mbox{$a_i(u_j) = 1$} if $i = j$ and is $0$ if $i\neq j.$ We write the dual of $u_1^{d_1}\ldots u_n^{d_n}$ as $a_1^{(d_1)}\ldots a_n^{(d_n)}$ where the parenthesized exponents are called \textit{divided powers}, and be careful that  in the corresponding situation over a field of characteristic $0$ in place of $\mathbb Z/2,$ $a_i^{(d_i)} = a_i^{d}/d!,$ which fits with the formula $a^{(d)}a^{(e)} = \binom{d+e}{d}a^{(d+e)}.$ This product gives a commutative graded algebra $\Gamma(a_1, \ldots, a_n)$ over $\mathbb Z/2$ called a \textit{divided power algebra}, where $\Gamma(a_1, \ldots, a_n) = H_*(V^{\oplus n}),$ and an element $a_1^{(d_1)}\ldots a_n^{(d_n)}$ in $H_*(V^{\oplus n})$ corresponding to a monomial $u_1^{d_1}\ldots u_n^{d_n}$ in $H^{*}(V^{\oplus n})$ is called \textit{$d$-monomial}.
%is generated by $a_1, \ldots, a_n,$ where $a_i$ is linear dual to $u_i\in H^{1}(V_n).$ 
Now, let $P_{\mathbb A}H_*(V^{\oplus n})$ be the subspace of $H_*(V^{\oplus n})$ consisting of all elements that are annihilated by all Steenrod squares of positive degrees. The general linear group $GL_n = GL(V^{\oplus n})$ acts  regularly on the classifying space $BV^{\oplus n}$ and therefore on $H^{*}(V^{\oplus n})$ and $H_*(V^{\oplus n}).$ This action commutes with that of the algebra $\mathbb A$ and so acts $\mathbb Z/2 \otimes _{\mathbb A} H^{*}(V^{\oplus n})$ and $P_{\mathbb A}H_*(V^{\oplus n}).$ For each $n\geq 0,$ Singer constructed in \cite{Singer} a linear transformation from $P_{\mathbb A}H_*(V^{\oplus n})$ to the $n$-th cohomology group ${\rm Ext}_{\mathbb A}^{n, n+*}(\mathbb Z/2, \mathbb Z/2)$ of $\mathbb A$, which commutes with two $Sq^{0}$'s on $P_{\mathbb A}H_*(V^{\oplus n})$ and ${\rm Ext}_{\mathbb A}^{n, n+*}(\mathbb Z/2, \mathbb Z/2)$ (see Boardman \cite{J.B} and Minami \cite{N.M2} for more about this). He shows that this map factors through the quotient of its domain's $GL_n$-coinvariants to give rise the so-called \textit{cohomological transfer of rank n}\\[1mm]
$\begin{array}{ll}
 \phi_n^{*}(\mathbb Z/2)&: \mathbb Z/2 \otimes _{GL_n}P_{\mathbb A}H_{*}(V^{\oplus n})\\
&\hspace{2.9cm}\longrightarrow {\rm Ext}_{\mathbb A}^{n, n+*}(\mathbb Z/2, \mathbb Z/2).
\end{array}$\\[1mm]
The domain of this transfer is dual to the space of $GL_n$-invariants $(\mathbb Z/2 \otimes _{\mathbb A} H^{*}(V^{\oplus n}))^{GL_n}.$ It is to be noted that $\phi_n^{*}(\mathbb Z/2)$ is induced over the $E_2$-term of the ASS by the geometrical transfer map $\Sigma^{\infty}(B(V^{\oplus n})_{+})\longrightarrow \Sigma^{\infty}(\mathbb S^{0})$ in stable homotopy theory (see also Mitchell \cite{S.M}). The work of Minami \cite{N.M2} indicated that these transfers play a key role in finding permanent cycles in the ASS. In the second cohomology groups of $\mathbb A$, following Mahowald \cite{Mahowald} and Lin-Mahowald \cite{L.M}, the classes $h_1h_j$ for $j\geq 3$ and $h_j^{2}$ for $0\leq j\leq 5,$ are known to be the permanent cycles in the ASS. In 2016, Hill, Hopkins, and Ravenel \cite{Hill} showed that when $j\geq 7,$ the class $h_j^{2}$ is not a permanent cycle in the ASS. It is surprising that so far there is no answer for $j = 6.$ The question of whether these $h_j^{2}$ are the permanent cycles in the ASS or not is called \textit{Kervaire invariant problem} in literature  \cite{W.B}. This is one of the oldest unresolved issues in Differential and Algebraic topology. 

Direct calculating the value of $\phi_n^{*}(\mathbb Z/2)$ on any non-zero element is a hard work. It has been demonstrated that $\phi_n^{*}(\mathbb Z/2)$ is an isomorphism for $n \leq 2$ by Singer himself \cite{Singer}, and $n = 3$ by Boardman \cite{J.B}. Most notably,   Singer sets up a hypothesis in the same paper \cite{Singer}  that \textit{$\phi_n^{*}(\mathbb Z/2)$ is a monomorphism}, but this is still not confirmed, for all cohomological degrees $n\geq 4.$ The cases $n = 4,\, 5$ are our concern in this paper. Besides Singer's transfer homomorphism, the lambda algebra $\Lambda$ of Bousfield et al. \cite{Bousfield} is also a relatively efficient tool to compute the cohomology of the Steenrod algebra. Recall that $\Lambda$ is the quotient of  the graded tensor algebra over $\mathbb Z/2$ on symbols $\lambda_i$ for $i\geq -1,$ modulo the two-sided ideal generated by $\lambda_s\lambda_k - \sum_{j}\binom{j-k-1}{2j-s}\lambda_{s+k-j}\lambda_j,$ for any $s,\, k\geq -1$ by the right ideal generated by $\lambda_{-1}.$ An interesting representation in the algebra $\Lambda$ of the algebraic transfer, established by Ch\ohorn n and H\`a \cite{C.H}, is a $\mathbb Z/2$-linear map $\psi_n$ from $P_{\mathbb A}H_{*}(V^{\oplus n})$ to a subspace of $\Lambda$ spanned by all monomials of length $n$ in all the monomials in $\lambda_i.$ The authors showed that the image of an element \mbox{$\zeta\in P_{\mathbb A}H_{*}(V^{\oplus n})$} under $\psi_n$ is a cycle in $\Lambda$ and $[\psi_n(\zeta)] = \phi_n^{*}(\mathbb Z/2)([\zeta]).$ Note also that this result is a dual version of the one in H\uhorn ng \cite{Hung}.

The Singer transfer we are discussing is closely related to the \textit{hit problem} in literature  \cite{F.P} of determination of a minimal generating set for the unstable $\mathbb A$-module $H^{*}(V^{\oplus n}).$ The reader can find an excellent list of publications about this problem in the works by Kameko \cite{M.K}, Mothebe-Uys \cite{M.M}, the present writer \cite{D.P2, D.P4, D.P7}, Singer \cite{Singer2}, Sum \cite{N.S1, N.S3}, Walker-Wood \cite{W.W}, Wood \cite{R.W} and others. Hit problems are motivated by several problems in Topology and Algebra. It was completely studied for the cases $n\leq 4$ by Peterson \cite{F.P}, Kameko's thesis \cite{M.K} and Sum \cite{N.S1}. Nevertheless, the general answer seems to be out of reach with the present techniques. Therefore, it is renowned as a difficult problem, even with the help of a computer. In fact, when $\mathbb Z/2$ is a trivial $\mathbb A$-module, solving the hit problem is equivalent to determining the "cohits" \mbox{$\mathbb Z/2 \otimes_{\mathbb A} H^{*}(V^{\oplus n})$} as a graded vector space, or more generally as a graded module over the group algebra $\mathbb Z/2[GL_n].$ Frank Peterson \cite{F.P} conjectured that $\mathbb Z/2 \otimes_{\mathbb A} H^{*}(V^{\oplus n}) = 0$ unless $\alpha(n+*)\leq n,$ where $\alpha(k)$ is the number of $1$'s in the dyadic expansion of a positive integer $k.$ His motivation for this was to prove that if $\mathcal M$ is a smooth manifold of dimension $*$ such that all products of length $n$ of Stiefel-Whiney classes of its nomal bundle vanish, then either $\alpha(*)\leq n$ or $\mathcal M$ is cobordant to zero. The conjecture was established by Wood \cite{R.W}. Therefrom,  to study $\mathbb Z/2 \otimes_{\mathbb A} H^{*}(V^{\oplus n})$ in each $n$ and degree \mbox{$* \geq 0$},  it suffices by Peterson's conjecture and iteration of  the Kameko map \cite{M.K} to consider degrees $*$ in the following "generic" form:
\begin{equation}\label{ct}
k(2^{t}-1) + \ell.2^{t},\ \mbox{ for $k,\, t, \, \ell \geq 0,\, \mu(\ell) < k \leq n$},
\end{equation}
where $\mu(\ell)$ is a smallest number $r\in \mathbb N$ such that $\alpha(\ell + r)\leq r.$ By Minami \cite{N.M}, hit problems are also considered as an useful tool for studying permanent cycles in the ASS.

In the present work, we explicitly determine the structure of the coinvariant $\mathbb Z/2 \otimes _{GL_n}P_{\mathbb A}H_{*}(V^{\oplus n})$ and the behavior of the cohomological transfer of ranks 4 and 5 in some generic degrees of the form \eqref{ct} by using techniques of the hit problem and the representation in the algebra $\Lambda$ of these transfers.

\section{Main results}
To begin with, we remark that by Sum \cite{N.S1}, it is enough to depict the behavior of the fourth transfer in the following degrees $d$:\\[1mm]
$ \begin{array}{lll}
{\rm (i)}  & d &= 2^{t+1} - r, \ \mbox{for $1\leq r\leq 3,$} \\
{\rm (ii)} & d &= 2^{t+s+1} +2^{t + 1}-3,\\
{\rm (iii)} & d &= 2^{t+s} +2^{t}-2,\\
{\rm (iv)} & d &= 2^{t+s+u} + 2^{t+s} + 2^{t} - 3,
\end{array}$\\[1mm]
whenever $r,\, t,\, s,$ and $u$ are positive integers. It is not difficult to check that the above degrees can be rewritten as \eqref{ct}. The cases of (i) are known by Sum \cite{N.S2}. The results for (ii) and (iii) have been partially probed in \cite{D.P5, D.P6}. This note is to investigate the case (iv).

Now, let us state the main results of this text. We first study the fourth cohomological transfer in degrees \mbox{$d_{t,\,s,\, u} := 2^{t+s+u} + 2^{t+s} + 2^{t} - 3$}. To make this, we give an explicit description of a minimal generating set for the domain of $\phi_4^{*}(\mathbb Z/2)$ in degree $d_{t,\,s,\, u},$ and obtain the following, which is proved in many steps by using computational techniques of the hit problem of four variables and some previous results by the present author \cite{D.P5, D.P6} and Sum \cite{N.S2}.
%The proofs of this theorem are made in many steps and are extremely technical. 
\begin{thm}\label{dlc1}
The domain of $\phi_4^{*}(\mathbb Z/2)$ in degree $d_{t,\,s,\, u}$ is determined by
$$ \begin{array}{ll}
\medskip
&\mathbb Z/2 \otimes _{GL_4}P_{\mathbb A}H_{d_{t,\,s,\, u}}(V^{\oplus 4}) \\
&=\left\{\begin{array}{ll}
0&\mbox{if $s = 1$, $u = 1$ and $t\geq 1$},\\
%\langle [\zeta_{1,\, 2,\, 1}] \rangle&\mbox{if $s = 2$, $u = 1$ and $t= 1$},\\
0 &\mbox{if $s = 2$, $u = 1$ and $t\geq 2$},\\
0&\mbox{if $s\geq 3$, $u=1$ and $t\geq 1$},\\
0&\mbox{if $s = 1$, $u = 2$ and $t= 1$},\\
\langle [\zeta_{t,\, 1,\, 2}] \rangle,&\mbox{if $s = 1$, $u = 2$ and $t\geq 2$},\\
0&\mbox{if $s = 1$, $u\geq 3$ and $t\geq 1$},\\
\langle [\zeta_{1,\, 2,\, u}] \rangle,&\mbox{if $s = 2$, $u\geq 1$ and $t =  1$},\\
0&\mbox{if $s \geq 3$, $u\geq 2$ and $t= 1$},\\
\langle [\zeta_{t,\, s,\, u}] \rangle,&\mbox{if $s \geq 2$, $u\geq 2$ and $t\geq 2$},\\
%0, &\mbox{otherwise},
\end{array}\right.
\end{array}$$
where\\[1mm] 
$\begin{array}{ll}
\zeta_{t,\, 1,\, 2}&= a_1^{(0)}a_2^{(2^{t+2}-1)}a_3^{(2^{t+2}-1)}a_4^{(3.2^{t}-1)}\\
&\quad + a_1^{(0)}a_2^{(2^{t+2}-1)}a_3^{(5.2^{t}-1)}a_4^{(2^{t+1}-1)}\\
&\quad + a_1^{(0)}a_2^{(6.2^{t}-1)}a_3^{(3.2^{t}-1)}a_4^{(2^{t+1}-1)} \\
\medskip
&\quad +  a_1^{(0)}a_2^{(7.2^{t}-1)}a_3^{(2^{t+1}-1)}a_4^{(2^{t+1}-1)},\\
\zeta_{1,\, 2,\, u} &=  a_1^{(2^{u+3}-1)}a_2^{(3)}a_3^{(3)}a_4^{(2)}\\
&\quad + a_1^{(2^{u+3}-1)}a_2^{(3)}a_3^{(4)}a_4^{(1)}\\
&\quad + a_1^{(2^{u+3}-1)}a_2^{(5)}a_3^{(2)}a_4^{(1)} \\
\medskip
&\quad+  a_1^{(2^{u+3}-1)}a_2^{(6)}a_3^{(1)}a_4^{(1)},\\
\zeta_{t,\, s,\, u} &= a_1^{(0)}a_2^{(2^{t}-1)}a_3^{(2^{s+t}-1)}a_4^{(2^{s+t+u}-1)}.
\end{array}$
\end{thm}
It would be interesting also to see that by this theorem, we have isomorphisms:\\[1mm]
$ \mathbb Z/2 \otimes _{GL_4}P_{\mathbb A}H_{d_{t,\,1,\, 2}}(V^{\oplus 4})\cong \mathbb Z/2,\ (t\geq 2),\\
 \mathbb Z/2 \otimes _{GL_4}P_{\mathbb A}H_{d_{1,\,2,\, u}}(V^{\oplus 4})\cong \mathbb Z/2,\ (u\geq 1),\\
 \mathbb Z/2 \otimes _{GL_4}P_{\mathbb A}H_{d_{t,\,s,\, u}}(V^{\oplus 4}) \cong \mathbb Z/2,\, (t\geq 2, s\geq 2, u\geq 2).$\\[1mm]
 Since the elements $\zeta_{t,\, 1,\, 2},$ $\zeta_{1,\, 2,\, u}$ and $\zeta_{t,\, s,\, u}$ belong to $P_{\mathbb A}H_{d_{t,\,1,\, 2}}(V^{\oplus 4}),$  $P_{\mathbb A}H_{d_{1,\,2,\, u}}(V^{\oplus 4})$ and $P_{\mathbb A}H_{d_{t,\,s,\, u}}(V^{\oplus 4}),$ respectively, $\psi_4(\zeta_{t,\, 1,\, 2}),$ $\psi_4(\zeta_{1,\, 2,\, u})$ and $\psi_4(\zeta_{t,\, s,\, u})$ are cycles in the lambda algebra $\Lambda$. Furthermore, using the representation in $\Lambda$ of the fourth transfer map, it may be concluded that\\[1mm]
$\begin{array}{ll}
&\phi_4^{*}(\mathbb Z/2)([\zeta_{t,\, 1,\, 2}]) = [\psi_4(\zeta_{t,\, 1,\, 2})]\\
&= [\lambda_0\lambda_{2^{t+2}-1}^{2}\lambda_{2^{t+1} + 2^{t}-1}]\\
&  =  h_0c_t\in {\rm Ext}_{\mathbb A}^{4, 4+d_{t,\,1,\, 2}}(\mathbb Z/2, \mathbb Z/2),\\
&\phi_4^{*}(\mathbb Z/2)([\zeta_{1,\, 2,\, u}]) = [\psi_4(\zeta_{1,\, 2,\, u})]\\
&= [\lambda_{2^{u+3}-1}\lambda_3^{2}\lambda_2] \\
& =  h_{u+3}c_0\in {\rm Ext}_{\mathbb A}^{4, 4+d_{1,\,2,\, u}}(\mathbb Z/2, \mathbb Z/2),\\
&\phi_4^{*}(\mathbb Z/2)([\zeta_{t,\, s,\, u}]) = [\psi_4(\zeta_{t,\, s,\, u})] \\
& = [\lambda_0\lambda_{2^{t}-1}\lambda_{2^{s+t}-1}\lambda_{2^{s+t+u}-1}]\\
&=  h_0h_th_{s+t}h_{s+t+u}\in {\rm Ext}_{\mathbb A}^{4, 4+d_{t,\,s,\, u}}(\mathbb Z/2, \mathbb Z/2).
\end{array}$\\[1mm]
These, together with Theorem \ref{dlc1} and a fact of the fourth cohomology groups ${\rm Ext}_{\mathbb A}^{4, 4+d_{t,\,s,\, u}}(\mathbb Z/2, \mathbb Z/2)$ (see Lin \cite{W.L}), yields the following.

\begin{corl}\label{hq1}
The algebraic transfer is an isomorphism in bidegree $(4, 4+d_{t,\,s,\, u})$ for all $t,\, s,\, u.$
\end{corl}

Next, we are going to survey the behavior of the transfer homomorphism of rank 5 in degrees of the form \eqref{ct} with $k = n = 5$ and $\ell = 50.$ To accomplish this, the domains of $\phi_5^{*}(\mathbb Z/2)$  at these degrees are computed as follows.

\begin{thm}\label{dlc2}
Let $d_{t} := 5(2^{t}-1) + 50.2^{t}$ with $t$ an arbitrary non-negative integer. Then, the coinvariant spaces $\mathbb Z/2\otimes_{GL_5} P_{\mathbb A}H_{d_t}(V^{\oplus 5})$ are trivial. 
\end{thm}

The proof of the theorem is quite long and complicated. We give a brief description for this: Firstly, we compute explicitly the monomial bases for the spaces of indecomposable elements $\mathbb Z/2\otimes_{\mathbb A}H^{*}(V^{\oplus 5})$ in degrees $d_t.$ An effective approach, based on the iterated Kameko squaring operation \cite{M.K} coupled with the results by Mothebe-Uys \cite{M.M}, T\'in \cite{Tin} and some linear transformations in Sum \cite{N.S1}, is applied in our computations Additionally, we also used the MAGMA computer algebra to verify the results. In particular, these calculations confirmed Sum's conjecture \cite{N.S3} for the relationship between $\mathbb A$-generators of the algebras $H^{*}(V^{\oplus 4})$ and $H^{*}(V^{\oplus 5})$ in the above degrees. Next, using those monomial bases, we may claim that the invariants  $(\mathbb Z/2\otimes_{\mathbb A}H^{*}(V^{\oplus 5}))^{GL_5}$ are trivial in degrees $d_t,$ from which the theorem follows from the fact that $\mathbb Z/2\otimes_{GL_5} P_{\mathbb A}H_{*}(V^{\oplus 5})$ is isomorphic to $(\mathbb Z/2\otimes_{\mathbb A}H^{*}(V^{\oplus 5}))^{GL_5}.$ 

Now according to \cite{W.L}, it is easily seen that\\[1mm]
${\rm Ext}_{\mathbb A}^{5, 5+d_{t}}(\mathbb Z/2, \mathbb Z/2) = \langle h_0h_{t+1}h_{t+2}h_{t+4}h_{t+5} \rangle = 0,$ for arbitrary $t\geq 0,$ and so, by Theorem \ref{dlc2}, we immediately obtain

\begin{corl}\label{hq2}
The Singer transfer is a trivial isomorphism in bidegree $(5, 5+d_{t})$ for every non-negative integer $t.$
\end{corl}

Thus, Corollaries \ref{hq1} and \ref{hq2} favor the Singer conjecture in bidegrees $(4, 4+d_{t,\,s,\, u})$ and $(5, 5+d_{t})$ for any $t,\, s,\, u.$

Detailed proofs of all the results of this note will be published elsewhere.

\section{Conclusion}

Although our work does not apparently lead to either a proof or a refutation of the Singer conjecture in general, we feel that it represents an interesting note about an application of hit problem and the lambda algebra for studying this conjecture. Perhaps a continuation of our methods may provefruitful.  It is our belief that further research canspring from these ideas.

\begin{acknow}
This work belongs to the Research Project 2022.  The author would like to express warm gratitude to Prof. Nguy\~\ecircumflex n Sum for many enlightening e-mail exchanges.
\end{acknow}

\end{document}